%% file: oneholefrompants.arXiv.tex
\documentclass[psamsfonts]{amsart}

\usepackage{amssymb, amscd, graphicx} 
\usepackage[all]{xy}

\include{oneholefrompants.arXiv.bbl}

\newcommand{\g}{\mathfrak{g}}
\newcommand{\C}{\mathbb{C}}
\newcommand{\tr}[1]{\mathrm{tr}\left(#1\right)}
\newcommand{\xb}{\mathbf{X}}

\newcommand{\wb}{\mathbf{W}}

\newcommand{\aq}{/\!\!/}
\newcommand{\X}{\mathfrak{X}}
\newcommand{\R}{\mathfrak{R}}
\newcommand{\hm}{\mathrm{Hom}}
\newcommand{\SL}{\mathrm{SL}(3,\C)}

\newcommand{\F}{\mathtt{F}}
\newcommand{\G}{\mathfrak{G}}
\newcommand{\xt}{\mathtt{x}}
\newcommand{\yt}{\mathtt{y}}

\newcommand{\bt}{\mathtt{b}}

\newcommand{\gt}{\mathtt{g}}

\newcommand{\ti}[1]{t_{(#1)}}

\newcommand{\Ad}{\mathrm{Ad}}

\newtheorem{theorem}{Theorem}

\newtheorem{corollary}[theorem]{Corollary}

\newtheorem{commnt}[theorem]{Comment}

\title[One-Holed Torus from Pants]{Obtaining the One-Holed Torus from Pants:
               Duality in an $\SL$-Character Variety }
\author[S. Lawton]{Sean Lawton}
\address{Mathematics Department,Instituto Superior T\'ecnico, Lisbon, Portugal}
\email{slawton@math.ist.utl.pt}
\urladdr{http://www.math.ist.utl.pt/$\sim$slawton}
\date{\today}

\begin{document}

\bibliographystyle{amsplain}

\maketitle

\begin{abstract}
The $\SL$-representation variety $\R$ of a free group $\F_r$ arises naturally by considering surface group representations
for a surface with boundary.  There is a $\SL$-action on the coordinate ring of $\R$.  The geometric points of the subring of invariants of this action is an affine
variety $\X$.  The points of $\X$ parametrize isomorphism classes of completely reducible representations.  The coordinate ring $\C[\X]$ is a complex Poisson algebra with respect to 
a presentation of $\F_r$ imposed by the surface.  In previous work, we have worked out the bracket on all generators when the surface is a three-holed sphere and when the surface is a one-holed torus.  In this paper, we show how the symplectic leaves corresponding to these two different Poisson structures on $\X$ relate to each other.  In particular, they are symplectically dual at a generic point.  Moreover, the topological gluing map which turns the three-holed sphere into the one-holed torus induces a rank preserving Poisson map on $\C[\X]$.
\end{abstract}


\section{Introduction}
In \cite{L3} we describe two competing Poisson structures on the variety of characters of representations of a rank 2 free group into $\SL$.  The purpose of this paper is to show that these two structures generically define symplectically dual symplectic leaves, and that a natural topological mapping relates the two Poisson structures in a non-trivial fashion.

For the remainder of this section we briefly describe character varieties and their smooth stratum's foliation by complex symplectic submanifolds.  In these terms we formulate our main theorems.  In Section 2, we describe in further detail past results necessary to make sense of the discuss at hand.  In particular, for the three-holed sphere and the one-holed torus we explicitly review the algebraic structure and the Poisson structure of the character variety.  Lastly, in Section 3 we restate and prove our main theorems.

\subsection{Algebraic Structure of $\X(\Sigma_{n,g})$}
Let $\Sigma_{n,g}$ be a compact, connected, oriented surface of genus $g$ with $n>0$ disks removed.  If $g=0$ we assume $n\geq 3$.
Its fundamental group has the following presentation:
 $$\pi_1(\Sigma_{n,g},*) =\{\xt_1,\yt_1,...,\xt_g,\yt_g,\bt_1,...,\bt_n\ | \ \Pi_{i=1}^{g}[\xt_i,\yt_i]\Pi_{j=1}^{n}\bt_j=1\}.$$ 
The group $\F_r:=\pi_1(\Sigma_{n,g},*)$ is always free of rank $r=2g+n-1$ since $\Sigma_{n,g}$ retracts to a wedge of $2g+n-1$ circles.  Let $\G=\SL$ and let $\{\gt_1,...,\gt_r\}$ be generators of $\F_r$.  The representation variety $\R=\hm(\F_r,\G)$ is bijectively equivalent to $\G^{\times r}$ given by evaluation:  $$\rho\mapsto
(\rho(\gt_1),...,\rho(\gt_r)),$$ and so inherits the structure of a smooth affine variety from $\G$.  The coordinate ring of $\G$ is the complex polynomial ring in $9$ indeterminates subject to the irreducible relation $\det(\xb)-1$ where $\xb=(x_{ij})$ is a \emph{generic matrix} and $x_{ij}$ are the $9$ indeterminates.
There is a polynomial action of $\G$ on the coordinate ring of $\R$, denoted by $\C[\R]$, by conjugation in $r$ \emph{generic matrices}; that is, for $g\in \G$
and $f\in \C[\R]$: $$g\cdot f(\xb_1,...,\xb_r)=f(g^{-1}\xb_1 g,..,g^{-1}\xb_r g).$$
Work of Procesi from 1976 \cite{P1} implies that the ring of invariants $\C[\R]^\G$ is generated by $\{\tr{\wb}\ |\ \mathtt{w}\in \F_r, \
|\mathtt{w}|\leq 6\}$.  Here $\wb$ is the word $\mathtt{w}$ in $\F_r$ with its letters replaced by generic matrices.
Thus,  $\C[\R]^\G$ is a finitely generated domain, and so its geometric points are an irreducible algebraic set, $\X(\Sigma_{n,g})=\mathrm{Spec}(\C[\R]^{\G})=\R\aq\G$, called the $\G$-\emph{character variety of} $\F_r$.  We note that the quotient notation just used means that it is a categorical quotient for the $\G$-action (see \cite{Do} or \cite{MFK}).

\subsection{The Boundary Map and Foliation of the Top Stratum}

The coordinate ring of $\G\aq\G$ is $$\C[\G\aq\G]=\C[\tr{\xb},\tr{\xb^{-1}}].$$  So $\G\aq\G=\C^2$ which we parametrize by coordinates $(\tau_{(1)}, \tau_{(-1)})$.  We then define the boundary map
$$\mathfrak{b}_i:\X=\R\aq\G=\hm(\pi_1(\Sigma_{n,g},*),\G)\aq\G\longrightarrow\G\aq\G$$ by sending a representation class, $[\rho]\mapsto
[\rho_{|_{\bt_i}}]=(\tau^i_{(1)}, \tau^i_{(-1)})$, to the class corresponding to the restriction of $\rho$ to the boundary $\bt_i$.
Subsequently we define $$\mathfrak{b}_{n,g}=(\mathfrak{b}_1,...,\mathfrak{b}_n):\X=\G^{\times r}\aq\G\longrightarrow (\G\aq\G)^{\times n}.$$
The map $\mathfrak{b}_{n,g}$ depends on the surface, not only its fundamental group.  We refer to it as a \emph{peripheral structure}, and the pair $(\X,\mathfrak{b}_{n,g})$ as the \emph{relative character variety}. 

Let $\mathfrak{L}=\mathfrak{b}_{n,g}^{-1}((\tau^1_{(1)}, \tau^1_{(-1)}),...,(\tau^n_{(1)}, \tau^n_{(-1)})),$ and $\mathcal{X}$ be the complement of the singular locus (a proper closed sub-variety) in
$\X$.  So $\mathcal{X}$ is a complex manifold that is dense in $\X$.  At regular values of $\mathfrak{b}_{n,g}$ (these are generic since $\mathfrak{b}_{n,g}$ is dominant), $\mathfrak{L}\cap\mathcal{X}$ is a submanifold of dimension $8r-8-2n=16(g-1)+6n$.  It is shown in \cite{L3} that the union of these \emph{leaves}, $\mathcal{L}=\mathfrak{L}\cap\mathcal{X}$, foliate $\mathcal{X}$ by complex symplectic submanifolds, making $\mathcal{X}$ a complex Poisson manifold.  This structure continuously extends over all of $\X$.

\subsection{The Quotient Map and Main Results}

There are two orientable surfaces with Euler characteristic $-1$; the three-holed sphere and the one-holed torus.  Both of these surfaces have fundamental groups free of rank 2.
Moreover, there is a natural topological quotient mapping (independent of orientation) which maps the three-holed sphere (hereafter referred to as pants) to the one-holed torus  $$q:\Sigma_{3,0}\to\Sigma_{1,1}.$$  In \cite{L3} we work out explicitly (with respect to a coordinate system for $\X$ and choices of orientation for the surfaces) the Poisson structures for the pants and the one-holed torus.

We now state our main theorems:

\begin{theorem}\label{mainone}Depending on the choice of orientation, $$q^*:\C[\X(\Sigma_{3,0})]\to \C[\X(\Sigma_{1,1})]$$ is generically a rank preserving Poisson (anti) morphism.\end{theorem}

Let $\mathcal{L}(\Sigma_{3,0})$ and $\mathcal{L}(\Sigma_{1,1})$ be generic symplectic leaves of $\X$.
 \begin{theorem}\label{maintwo}
$\mathcal{L}(\Sigma_{3,0})$ and $\mathcal{L}(\Sigma_{1,1})$ are (generically) transverse and so are symplectically dual to each other.
 \end{theorem}

\subsection*{Acknowledgments}
I would like to thank my adviser Bill Goldman for introducing me to this subject.  Also, this work was inspired by a visit to McMaster's University while I was visiting Hans Boden.  I would like to thank him and McMaster's for that opportunity.  Additionally, I would like to thank Instituto Superior T\'ecnico for hosting me while I was working on the draft of this paper.

\section{Past Results and Background}

In this section we very briefly review some of the results from \cite{L2, L3} that we will need to prove our theorems.

\subsection{Symplectic and Poisson Structure on $\X(\Sigma_{n,g})$}
In 1997, Guruprasad, Huebschmann, Jeffrey, and Weinstein \cite{GHJW} showed $\omega$ (in the following commutative diagram) defines a symplectic form on the leaf $\mathcal{L}$ defined in the introduction:

$$
\xymatrix{
H^1(S,\partial S;\g_\Ad) \times H^1(S;\g_\Ad) \ar[r]^-\cup   & H^2(S,\partial S;\g_\Ad\otimes \g_\Ad) \ar[d]^{\mathrm{tr}_*}\\
    &  H^2(S,\partial S;\C) \ar[d]^{\cap [Z]}\\
H^1_{\mathrm{par}}(S;\g_\Ad)\times H^1_{\mathrm{par}}(S;\g_\Ad)\ar[uu] \ar[r]^-\omega  & H_0(S;\C)=\C }
$$

With respect to this 2-form, we show in \cite{L3} Goldman's proof \cite{G8,G7} of the Poisson bracket (a Lie bracket and derivation) generalizes directly to relative cohomology.

Let $\Sigma$ be an oriented surface with boundary, and $\alpha, \beta \in \pi_1(\Sigma,*)$.  Let $\alpha\cap\beta$ be the set of (transverse) double point intersections of $\alpha$ and $\beta$.  Let $\epsilon(p,\alpha, \beta)$ be the oriented intersection number at $p\in \alpha\cap\beta$ and let $\alpha_p\in \pi_1(\Sigma,p)$ be the curve $\alpha$ based at $p$. 

In these terms the bracket is defined on $\mathbb{C}[\mathfrak{X}]$ by:

\begin{align*}
\{\mathrm{tr}(\rho(\alpha)),\mathrm{tr}(\rho(\beta))\}&=\sum_{p\in
\alpha\cap\beta}\epsilon(p,\alpha,\beta)\big(\mathrm{tr}(\rho(\alpha_p\beta_p))-\frac{1}{3}\mathrm{tr}(\rho(\alpha))\mathrm{tr}(\rho(\beta))\big).
\end{align*}

\subsection{Algebraic Structure of $\X(\Sigma_{3,0})$ and $\X(\Sigma_{1,1})$}
We now review the algebraic structure of $\C[\X]$ for the pants and one-holed torus and the corresponding Poisson structures in those cases.  Details are available in \cite{L2, L3}.

Let $$\C[t_{(1)},t_{(-1)},t_{(2)},t_{(-2)},t_{(3)},t_{(-3)},t_{(4)},t_{(-4)},
t_{(5)}]$$ be a freely generated complex polynomial ring,  and let
$$R=\C[t_{(1)},t_{(-1)},t_{(2)},t_{(-2)},t_{(3)},t_{(-3)},t_{(4)},t_{(-4)}].$$

Define the following ring homomorphism,
$$R[t_{(5)}]\stackrel{ \Pi}{\longrightarrow} \C[\X]$$ by
\begin{center}

\begin{tabular}{ll}
$t_{(1)}\mapsto\tr{\xb_1}$ & $t_{(-1)}\mapsto\tr{\xb_1^{-1}}$\\
$t_{(2)}\mapsto\tr{\xb_2}$& $t_{(-2)}\mapsto\tr{\xb_2^{-1}}$\\
$t_{(3)}\mapsto\tr{\xb_1\xb_2}$& $t_{(-3)}\mapsto\tr{\xb_1^{-1}\xb_2^{-1}}$\\
$t_{(4)}\mapsto\tr{\xb_1\xb_2^{-1}}$& $t_{(-4)}\mapsto\tr{\xb_1^{-1}\xb_2}$\\
$t_{(5)}\mapsto\tr{\xb_1\xb_2 \xb_1^{-1}\xb_2^{-1}}$
\end{tabular}
\end{center}

It can be shown using trace equations that $\Pi$ is surjective, and hence $$R[t_{(5)}]/\ker(\Pi)\cong \C[\X].$$  The Krull dimension of $\X$ is $8$ since generic orbits are $8$ dimensional.  Hence, $\ker(\Pi)$ is non-zero and principal.

Let $\mathbb{S}$ be the formal sum of the elements in the group generated by the permutations (in cycle notation) $$(1,2)(-1,-2)(4,-4) \text{ and }
(1,-1)(3,-4)(-3,4)$$ acting on the indices of the generators of $R[t_{(5)}]/\ker(\Pi)$.
The action is induced by the following elements of the $\mathrm{Out}(\F_2)$:
 $$\begin{array}{cc}
\mathfrak{t}=\left\{
\begin{array}{l}
\xt_1\mapsto \xt_2\\
\xt_2\mapsto \xt_1
\end{array}\right. &
\mathfrak{i}_1=\left\{
\begin{array}{l}
\xt_1\mapsto \xt_1^{-1}\\
\xt_2\mapsto \xt_2
\end{array}\right.
\end{array}$$

The group generated has order $8$ and is isomorphic to the dihedral group $D_4$.  In \cite{L2} we show

\begin{theorem}\label{coordinatering}
\begin{enumerate}
\item[{}]
\item  $\X=\G^{\times 2}\aq\G$ is a degree $6$ hyper-surface in $\C^9$
\item $\ker(\Pi)=(\ti{5}^2-P\ti{5}+Q)$ where $P,Q\in R$
\item There is a $D_4$-equivariant surjection (submersion) $m:\X\to \C^8$, generically $2$-to-$1$.
\item $P$ and $Q$ are given by:
\begin{align*}
P=&\mathbb{S}\Big(\frac{1}{8}\big(t_{(1)}t_{(-1)}t_{(2)}t_{(-2)}-4t_{(1)}t_{(-2)}t_{(-4)}+2t_{(1)}t_{(-1)}+2t_{(3)}t_{(-3)}\big)\Big)-3\\
Q=&\mathbb{S}\Big( \frac{1}{8}\big(2t_{(-2)}t_{(-1)}^2 t_{(1)}^2t_{(2)}+4t_{(1)}^2t_{(2)}^2t_{(3)} -4t_{(1)}^3t_{(-2)}t_{(2)}-\\
&8t_{(-4)}t_{(-2)}t_{(-1)}t_{(1)}^2-4t_{(4)}t_{(3)}t_{(2)}t_{(1)}t_{(-2)} +8t_{(1)}t_{(3)}t_{(-4)}^2 +\\
&8t_{(-4)}t_{(1)}t_{(2)}^2 -8t_{(3)}^2 t_{(2)}t_{(1)}+4t_{(4)}t_{(-3)}t_{(2)}^2 +t_{(-2)}t_{(-1)}t_{(2)}t_{(1)} +\\
&t_{(-3)}t_{(-4)}t_{(3)}t_{(4)}+4t_{(-3)}t_{(-1)}t_{(3)}t_{(1)}+4t_{(1)}^3 +4t_{(3)}^3+\\
&12t_{(-4)}t_{(-2)}t_{(1)}-12t_{(-4)}t_{(2)}t_{(3)} -12t_{(1)}t_{(-1)} -12t_{(3)}t_{(-3)}\big) \Big)+9 \end{align*} 
\end{enumerate}
\end{theorem}

\subsection{Poisson Structures of $\X(\Sigma_{3,0})$ and $\X(\Sigma_{1,1})$}\label{poissonstructure}
For our purposes a Poisson variety is an affine variety $\X$ over $\C$ endowed with a Lie bracket $\{,\}$ on its coordinate ring $\C[\X]$ that acts as a formal derivation (satisfies the Leibniz rule).  On the smooth strata of $\X$ (denoted by $\mathcal{X})$, it makes $\mathcal{X}$ a complex Poisson manifold in the usual sense (by Stone--Weierstrass).  For any such Poisson bracket there exists an exterior bi-vector field $\mathfrak{a}\in \Lambda^{2}(T\mathcal{X})$ whose restriction to symplectic leaves is given by the symplectic form as $\{f,g\}=\omega(H_g,H_f)$ (where $H_f=\{f,\cdot \}$ is called the Hamiltonian vector field).  Let $f,g\in \C[\X]$.  Then with respect to interior multiplication $\{f,g\}=\mathfrak{a}\cdot df\otimes dg.$  In local coordinates $(z_1,...,z_k)$ it takes the form 
$$\mathfrak{a}=\sum_{i,j}\mathfrak{a}_{i,j}\frac{\partial}{\partial z_i}\land 
\frac{\partial}{\partial z_j}$$ and so 
\begin{align*}
\{f,g\}&=\sum_{i,j}\left(\mathfrak{a}_{i,j}\frac{\partial}{\partial z_i}\land \frac{\partial}{\partial z_j}\right)\cdot\left(\frac{\partial f}{\partial z_i} dz_i \otimes \frac{\partial g}{\partial 
z_j}dz_j\right)\\
&=\sum_{i,j}\mathfrak{a}_{i,j}\left(\frac{\partial f}{\partial z_i}\frac{\partial g}{\partial z_j}-
\frac{\partial f}{\partial z_j}\frac{\partial g}{\partial z_i}\right).
\end{align*} 

Denote the bi-vector associated to $\X(\Sigma_{n,g})$ by $\mathfrak{a}(\Sigma_{n,g})$.  In the case of the pants or the one-holed torus, there are 9 generating functions of $\C[\X]$:  $\ti{\pm i}$ for $1\leq i \leq 4$ and $\ti{5}$.  Since the bi-vector is a Lie bracket and a derivation, its formulation is in these terms.  Let $\mathfrak{a}_{i,j}=\{\ti{i},\ti{j}\}$.  In \cite{L3} we show the following two structure theorems.

\begin{theorem}
The Poisson bi-vector on $\X(\Sigma_{3,0})$ is 
$$\mathfrak{a}(\Sigma_{3,0})=(P-2\ti{5})\frac{\partial}{\partial \ti{4}}\land \frac{\partial}{\partial \ti{-4}} +(1-\mathfrak{i})\left(\mathfrak{a}_{4,5}\frac{\partial}{\partial \ti{4}}\land \frac{\partial}{\partial \ti{5}} \right),$$ where

$\mathfrak{a}_{4,5}=\frac{\partial}{\partial \ti{-4}}(Q-\ti{5}P)$
and where $\mathfrak{i}=\mathfrak{i}_1\mathfrak{t}\mathfrak{i}_1\mathfrak{t}$ is the mapping $\xt_i\mapsto \xt_i^{-1}$.

\end{theorem}

In $D_4$ define $\mathfrak{i}_2=\mathfrak{t}\mathfrak{i}_1\mathfrak{t}$; the mapping which sends $\xt_2\mapsto \xt_2^{-1}$.

Additionally, define the following elements of the group ring of $D_4$:
\begin{itemize}
\item $\Sigma_{1}=1+\mathfrak{i}-\mathfrak{i}_1-\mathfrak{i}_2$
\item $\Sigma_{2}=1+\mathfrak{i}-\mathfrak{t}-\mathfrak{i}\mathfrak{t}$.
\end{itemize}

Note $$\frac{1}{2}\Sigma_1\Sigma_2=1+\mathfrak{i}-\mathfrak{i}_1-\mathfrak{i}_2-\mathfrak{t}-\mathfrak{i}\mathfrak{t}+\mathfrak{i}_1\mathfrak{t}+\mathfrak{i}_2\mathfrak{t}.$$

Then after doing 28 calculations and observing symmetry along the way, we conclude:

\begin{theorem}
The Poisson bi-vector field on $\X(\Sigma_{1,1})$ is 
\begin{align*}\mathfrak{a}(\Sigma_{1,1})=&\Sigma_{1}\bigg(\mathfrak{a}_{1,2}\frac{\partial}{\partial \ti{1}}\land \frac{\partial}{\partial \ti{2}}\bigg)
+\Sigma_2\left(\mathfrak{a}_{3,4}\frac{\partial}{\partial \ti{3}}\land \frac{\partial}{\partial \ti{4}}\right)\\
&+\frac{1}{2}\Sigma_1\Sigma_{2}\bigg(\mathfrak{a}_{1,3}\frac{\partial}{\partial \ti{1}}\land \frac{\partial}{\partial \ti{3}}
+\mathfrak{a}_{1,-3}\frac{\partial}{\partial \ti{1}}\land \frac{\partial}{\partial \ti{-3}}\bigg),\end{align*}

where: 
\begin{itemize}
\item $\mathfrak{a}_{1,2}=\ti{3}-\frac{1}{3}\ti{1}\ti{2}$
\item $\mathfrak{a}_{1,3}=\frac{2}{3}\ti{1}\ti{3}-\ti{-1}\ti{2}+\ti{-4}$
\item $\mathfrak{a}_{1,-3}=-\ti{-2}+\frac{1}{3}\ti{1}\ti{-3}$
\item $\mathfrak{a}_{3,4}=-\ti{1}^2+\ti{-1}-\ti{-4}\ti{-2}-\ti{2}\ti{-3}+\ti{-1}\ti{2}\ti{-2}-\frac{1}{3}\ti{3}\ti{4}$.
\end{itemize}

\end{theorem}
\begin{commnt}
We note here that the orientations chosen on these surfaces are opposite.  Our presentation of the pants has the boundary on the outside whereas the one-holed torus has the same boundary (after projection) on the inside.  Since the orientations of the boundaries are the same the surfaces are ``inside-out'' with respect to each other, and so the orientations are reversed.  Consequently, if the quotient mapping taking the pants to the one-holed torus is to preserve orientations one of the above two bi-vectors must be multiplied by $-1$.
\end{commnt}

\section{Obtaining the Torus from Pants}
Let $q: \Sigma_{3,0}\to \Sigma_{1,1}$ be the quotient map given by identifying two of the boundaries (call them $b_1$ and $b_2$).  Let $x_0$ be a fixed base point.  And let $x_1\in b_1$ and $x_2\in b_2$ be also fixed.

Then the third boundary $b_3$ in $\Sigma_{3,0}$ is homotopic to $(b_1b_2)^{-1}$ in $\pi_1(\Sigma_{3,0},x_0)$.

Let $\gamma_1$ and $\gamma_2$ be paths from $x_0$ to $x_1$ and $x_0$ to $x_2$ respectively.  Then $q(\gamma_1\gamma_2^{-1}):=\beta$ is a non-trivial based loop in $\Sigma_{1,1}$.

Moreover, $(\gamma_1\gamma_2^{-1})b_1(\gamma_1\gamma_2^{-1})^{-1}$ is homotopic to $b_2^{-1}$ since $\gamma_1^{-1}b_1\gamma_1$ is homotopic to $\gamma_2^{-1}b_2^{-1}\gamma_2$ in $\Sigma_{1,1}$.

Therefore,  $$q_{\sharp}: \pi_1(\Sigma_{3,0},x_0)=\langle b_1,b_2,b_3\ | \ b_1b_2b_3=1\rangle \to \pi_1(\Sigma_{1,1},x_0)=\langle \alpha, \beta, \gamma \ | \ [\alpha,\beta]\gamma=1\rangle$$ is injective and given by $$b_1\mapsto \alpha,\ b_2\mapsto \beta\alpha^{-1}\beta^{-1}, \text{ and } b_3\mapsto [\alpha,\beta]^{-1}.$$

Consequently, $q_*:\X(\Sigma_{1,1})\to \X(\Sigma_{3,0})$ is given by $$[(A,B)]\mapsto [(A,BA^{-1}B^{-1})],$$  
and $q^*:\C[\X(\Sigma_{3,0})]\to \C[\X(\Sigma_{1,1})]$ is given by $f\mapsto f\circ q_*$.

To be concrete we write the assignments which determine $q^*$: $$\ti{1}\mapsto \ti{1},\ \ti{-1}\mapsto \ti{-1},\ \ti{2}\mapsto \ti{-1},\ \ti{-2}\mapsto \ti{1},$$
$$\ti{3}\mapsto \ti{5},\ \ti{-3}\mapsto \tr{A^{-1}BAB^{-1}}=P-\ti{5},$$
$$\ti{4}\mapsto \tr{ABAB^{-1}}= \ti{-4}\ti{-2}+\ti{-1}+\ti{-3}\ti{2}-\ti{-1}\ti{2}\ti{-2}+\ti{3}\ti{4},$$
 $$\ti{-4}\mapsto \tr{A^{-1}BA^{-1}B^{-1}}=\ti{-4}\ti{-3}+\ti{1}+\ti{3}\ti{-2}-\ti{1}\ti{2}\ti{-2}+\ti{2}\ti{4},$$

and $\ti{5}\mapsto \tr{ABA^{-1}B^{-1}A^{-1}BAB^{-1}}=$
\begin{align*}
&\ti{-3}^3-2 \ti{-2} \ti{-1} \ti{-3}^2+\ti{-4} \ti{2} \ti{-3}^2+\ti{1} \ti{4}\ti{-3}^2\\
&+\ti{-2}^2 \ti{-1}^2 \ti{-3}-\ti{-2} \ti{1} \ti{2}^2 \ti{-3}+\ti{-1}\ti{4}^2 \ti{-3}\\
&+\ti{-4}^2 \ti{-2} \ti{-3}-\ti{-4} \ti{-1} \ti{-3}+\ti{-2}^2 \ti{1}\ti{-3}\\
&-\ti{-4} \ti{-2} \ti{-1} \ti{2} \ti{-3}+2 \ti{1} \ti{2} \ti{-3}+\ti{-2} \ti{2}\ti{3} \ti{-3}\\
&-3 \ti{3} \ti{-3}+\ti{2}^2 \ti{4} \ti{-3}-3 \ti{-2} \ti{4}\ti{-3}\\
&-\ti{-2} \ti{-1} \ti{1} \ti{4} \ti{-3}+\ti{-4} \ti{3} \ti{4}\ti{-3}\\
&+\ti{-2}^3+\ti{4}^3+\ti{-2}^2 \ti{-1} \ti{1} \ti{2}^2-2 \ti{-2} \ti{1}\ti{4}^2\\
&+\ti{2} \ti{3} \ti{4}^2+2 \ti{-4} \ti{-2} \ti{1}-\ti{-2}^3 \ti{-1} \ti{1}\\
&-3\ti{-2} \ti{2}-\ti{-4} \ti{-2}^2 \ti{1} \ti{2}-\ti{-2} \ti{-1} \ti{1} \ti{2}\\
&+\ti{-4}\ti{-2}^2 \ti{3}+2 \ti{-2} \ti{-1} \ti{3}-\ti{-2}^2 \ti{-1} \ti{2} \ti{3}\\
&+\ti{-2}^2\ti{1}^2 \ti{4}-\ti{-2} \ti{-1} \ti{2}^2 \ti{4}+\ti{-2} \ti{3}^2 \ti{4}\\
&-3 \ti{-4}\ti{4}+\ti{-2}^2 \ti{-1} \ti{4}+\ti{-4} \ti{-2} \ti{2} \ti{4}\\
&+2 \ti{-1} \ti{2}\ti{4}-\ti{1} \ti{3} \ti{4}-\ti{-2} \ti{1} \ti{2} \ti{3} \ti{4}+3.
\end{align*}

This last three identities are a consequence of recursive trace reduction formulas (see \cite{L2}, \cite{L4}).

Let $\{,\}_1$ be the bracket corresponding to $\Sigma_{1,1}$ and let $\{,\}_3$ be the bracket corresponding to $\Sigma_{3,0}$.  Let $\mathcal{Q}$ be the image of $q^*$.
We now prove Theorem \ref{mainone} from the Introduction.
\begin{theorem}
\begin{enumerate}
\item[]
\item $\mathcal{Q}$ is a Poisson subalgebra of $\C[\X(\Sigma_{1,1})]$
\item $q^*$ is an anti-Poisson morphism; that is, $\{q^*(f),q^*(g)\}_1=-q^*\big(\{f,g\}_3\big)$
\item at a generic point of $\X$, $\mathrm{rank}\left(\mathfrak{a}(\Sigma_{1,1})\bigg|_{\mathcal{Q}}\right)=\mathrm{rank}\left(\mathfrak{a}(\Sigma_{3,0})\right)$
\end{enumerate}

\end{theorem}
\begin{proof}

First we note that (2) implies (1).

To prove (2), since $q^*$ is an algebra morphism and the bracket is a derivation, it is enough to verify it on all generators of the algebra.  One can use the explicit form of the mapping $q^*$ and the explicit form of the bi-vectors to verify the result.  However, since $q$ preserves transversality of cycles, double points, and does not affect orientation it follows that for any two cycles $\alpha$ and $\beta$ in $\Sigma_{3,0}$ used in computing the bi-vector $\mathfrak{a}(\Sigma_{3,0})$ we have:
\begin{align} \label{qbracket}
&q^*\big(\{\mathrm{tr}(\rho(\alpha)),\mathrm{tr}(\rho(\beta))\}_3\big)=\\
&\sum_{q(p)\in q(\alpha)\cap q(\beta)}\epsilon(q(p),q(\alpha),q(\beta))\big(\mathrm{tr}(\rho(q(\alpha)_{q(p)}q(\beta)_{q(p)}))-\frac{1}{3}\mathrm{tr}(\rho(q(\alpha)))\mathrm{tr}(\rho(q(\beta)))\big).\nonumber
\end{align}

However, as already noted earlier in Section \ref{poissonstructure}, the intersection numbers $\epsilon(q(p),q(\alpha),q(\beta))$ and $\epsilon(p,\alpha,\beta)$ must be reversed since the bracket computations carried out in \cite{L3} are with respect to opposite orientations on the surfaces $\Sigma_{3,0}$ and $\Sigma_{1,1}$.

Thus Equation $\eqref{qbracket}$ is exactly $-\{\tr{q(\alpha)},\tr{q(\beta)}\}_1$, as was to be shown.

To prove (3) we first note that the rank of a bi-vector is the rank of the anti-symmetric matrix of functions $(\mathfrak{a}_{ij})$.  Then, from (2), there are only three non-zero coefficients to $\mathfrak{a}(\Sigma_{1,1})$ after restricting to the image of $q^*$.  Namely, $\ti{5}$ and $P-\ti{5}$ are Casimirs for $\{,\}_1$, and since the mapping is Poisson and $\ti{\pm 1}$ are fixed, it follows that they are Casimirs in the images since they are Casimirs in the pre-image.  Thus we are left with the image generators $q^*(\ti{j})$ for $j=4,-4,5$.  Since the bi-vector on the Poisson subalgebra is exactly the induced one, we explicitly formulate the bi-vector matrix and compute its rank; finding it generically 2.  In particular, the matrix has the form:
$$\left(
\begin{array}{lll}
 0 & a & b \\
 -a & 0 & c \\
 -b & -c & 0
\end{array}
\right)$$
where $a=\{q^*(\ti{4}),q^*(\ti{-4})\}_1$, $b=\{q^*(\ti{4}),q^*(\ti{5})\}_1$, and $c=\{q^*(\ti{-4}),q^*(\ti{5})\}_1$.

However any matrix of this form has rank 2 as long as all three of $a,b,$ and $c$ are not 0, in which case the rank is 0.  By direct calculation one sees that all polynomials $a,b,c$ are in terms of \emph{only} algebraically independent generators, and so none of $a,b,c$ are generically 0 on $\X$.  So generically the rank is 2, and the rank of $\{,\}_3$ is two since the rank is equal to the dimension of a symplectic leaf.  Hence, the mapping is rank preserving as long as it is not completely degenerate (which is generically the case).
 \end{proof}

\begin{commnt}
Equation $\eqref{qbracket}$ used in the above argument shows much more.  For any two surfaces $\Sigma_{n_1,g_1}$ and $\Sigma_{n_2,g_2}$ with $n_1>n_2>0$ and $\chi(\Sigma_{n_1,g_1})=\chi(\Sigma_{n_2,g_2})$, there is a quotient mapping (identifying pairs of boundary components) $q:\Sigma_{n_1,g_1}\to\Sigma_{n_2,g_2}$ which gives an injection on fundamental groups and therefore gives a mapping of coordinate rings $q^*:\C[\X(\Sigma_{n_1,g_1})]\to\C[\X(\Sigma_{n_2,g_2})]$. This is true not only for $\SL$ but for any complex algebraic reductive Lie group $\G$.  The above argument (Equation $\eqref{qbracket}$) shows that if the orientations of the surfaces correspond to each other then $q^*$ is a Poisson mapping and if the orientations are opposite then it is a anti-Poisson morphism.  Consequently, the image of $q^*$ is a Poisson subalgebra of the codomain of $q^*$ in general.

It does not seem clear whether rank preserving is a general property or not.
\end{commnt}



Let $\mathcal{L}(\Sigma_{n,g})$ be a generic symplectic leaf of $\X(\Sigma_{n,g})$.  We now prove Theorem \ref{maintwo} from the Introduction.
 \begin{theorem}
$\mathcal{L}(\Sigma_{3,0})$ and $\mathcal{L}(\Sigma_{1,1})$ are transverse (generically).
 \end{theorem}

\begin{proof}

The mapping from Theorem \ref{coordinatering} $m:\X\to \C^8$ is given by $$(\ti{1},\ti{-1},...,\ti{4},\ti{-4},\ti{5})\mapsto (\ti{1},\ti{-1},...,\ti{4},\ti{-4}).$$  It is surjective, and since the first eight generators are algebraically independent, it is submersive as well.

This immediately implies that the mapping $\mathfrak{b}_{3,0}:\X(\Sigma_{3,0})\to \C^6$ given by $$(\ti{1},\ti{-1},...,\ti{4},\ti{-4},\ti{5})\mapsto (\ti{1},\ti{-1},...,\ti{3},\ti{-3})$$ is likewise surjective and submersive.  Consequently, for any $x\in \mathcal{X}$, $\mathcal{L}(\Sigma_{3,0})=\mathfrak{b}_{3,0}^{-1}(\mathfrak{b}_{3,0}(x))$ has dimension 2.  Moreover, we can locally parametrize this leaf by the coordinates $(\ti{4},\ti{-4})$ since the other six coordinates $\ti{\pm i}$ for $i=1,2,3$ are held constant and $\ti{5}$ is then determined by the defining relation $\ti{5}^2-P\ti{5}+Q$.  In particular, flows through these coordinates determine a dimension 2 subspace $T_*\mathcal{L}(\Sigma_{3,0})\subset T_*\mathcal{X}$ of the tangent space.

Now consider the mapping $\mathfrak{b}_{1,1}:\X(\Sigma_{1,1})\to \C^2$ given by $$(\ti{1},\ti{-1},...,\ti{4},\ti{-4},\ti{5})\mapsto (\ti{5},P-\ti{5}).$$  Note that this is in fact the correct mapping since $P=\tr{[A,B]}+\tr{[B,A]}$ (see \cite{L2}).  It is shown in \cite{L3} that the boundary mapping is always surjective if $g>0$.

However, $\mathfrak{b}_{1,1}$ may not be everywhere submersive; in particular, when $P^2-4Q=0$.  However, there is an open dense set of $\X$ where $d\mathfrak{b}_{1,1}$ is onto (call it $\mathcal{U}$), since $\mathfrak{b}_{1,1}$ is surjective and regular (see \cite{L3}).  We may assume $\mathcal{U}\subset \mathcal{X}$.

Take any $u\in \mathcal{U}$.  Then (as is shown in \cite{L3}) the leaves $\mathcal{L}_1:=\mathfrak{b}_{1,1}^{-1}(\mathfrak{b}_{1,1}(u))\cap \mathcal{X}$ and $\mathcal{L}_3:=\mathfrak{b}_{3,0}^{-1}(\mathfrak{b}_{3,0}(u))\cap \mathcal{X}$ are complex symplectic manifolds of dimensions 6 and 2 respectively.  Consequently, these leaves are \emph{properly} transverse; that is, $\mathrm{dim}\mathcal{L}_1+\mathrm{dim}\mathcal{L}_3=\mathrm{dim}\mathcal{X}$.

We now show $\mathrm{dim} \mathcal{L}_1\cap \mathcal{L}_3=0$ and $\mathcal{L}_1\cap \mathcal{L}_3\not=\emptyset$.  At an intersection point, $P=\ti{5}+\ti{-5}:=C$ and $Q=\ti{5}\ti{-5}:=D$ and $\ti{\pm 1},\ti{\pm 2}, \ti{\pm 3}$ are all fixed.  Moreover, solving $P=C$ for $\ti{4}$ generically gives \begin{align}\label{t4sub} \ti{4}\mapsto& \frac{1}{\ti{-4}-\ti{-1} \ti{2}}\bigg(C+\ti{-4} \ti{-2} \ti{1}-\ti{-1} \ti{1}-\nonumber \\ &\ti{-2} \ti{2}+\ti{-3} \ti{1} \ti{2}-\ti{-2} \ti{-1} \ti{1}\ti{2}-\ti{-3} \ti{3}+\ti{-2} \ti{-1} \ti{3}+3\bigg).\end{align}

Now, substituting this into $Q-D=0$ gives a monic degree six polynomial in the variable $\ti{-4}$ since the degree in $\ti{4}$ of $Q$ is 3.  So the intersection is non-empty and of dimension 0 (at most 6 discrete points) if $\ti{-4}-\ti{-1} \ti{2}\not=0$.  Otherwise, setting $\ti{-4}=\ti{-1} \ti{2}$ and substituting this into $Q-D=0$ gives a monic degree 3 polynomial in $\ti{4}$.  Either way, the intersection is non-empty and of dimension 0.

We claim that the tangent space to $\mathcal{L}_3$ locally can be determined by flows through $\{\ti{1},\ti{-1},\ti{2},\ti{-2},\ti{3},\ti{-3}\}$ by solving for $\ti{4}$ and $\ti{-4}$ in terms of $\ti{\pm 1},...,\ti{\pm 3}$ on an open subset since both $P$ and $Q$ are constant on $\mathcal{L}_3$.

Explicitly, substituting Equation \eqref{t4sub} into $Q-D=0$, where $\ti{\pm i}$ for $1\leq i\leq 3$ are now not fixed, locally and generically give $\ti{-4}$, and subsequently $\ti{4}$, as functions of $\ti{\pm i}$ for $1\leq i\leq 3$.  Thus the flows through $\ti{\pm i}$, $1\leq i\leq 3$ determine a full dimensional tangent space to $\mathcal{L}_3$ whenever $\ti{-4}-\ti{-1}\ti{2}\not=0$.

Switching the roles of $\ti{4}$ and $\ti{-4}$ give a like result for at any point where $\ti{4}-\ti{1} \ti{-2}\not=0$.

However, the tangent space to any point in $\mathcal{L}_1$ is given by the kernel to the mapping $M:=\big(\partial f_i/\partial \ti{\pm j}\big)$ where $f_1=\ti{5}^2-P\ti{5}+Q, f_2=\ti{5}-a,$ and $f_3=P-\ti{5}-b$, and $C=a+b$ and $D=ab$ (see \cite{Ha}).  This follows since these three functions define the leaf as an algebraic set cut out of $\C^9$.  At any smooth point in the leaf the dimension of the kernel is 6.  So whenever $\ti{4}-\ti{1} \ti{-2}\not=0$ and $\ti{-4}-\ti{-1} \ti{2}\not=0$ using $P=C$, $Q=D$, and $\ti{5}=a$ from above, this matrix has its entries rational functions of $\ti{\pm 1},..,\ti{\pm 3}$ alone.  On the other hand, if any smooth point also satisfies $\ti{4}-\ti{1} \ti{-2}=0$ and $\ti{-4}-\ti{-1} \ti{2}=0$ solving for $\ti{\pm 4}$ again gives $M$ as a matrix in these six variables.  Thus the flows through these six coordinate functions always determine the tangent space at a smooth point of $\mathcal{L}_1$.

Consequently, the span of the flows through $\{\ti{4},\ti{-4}\}$ and $\{\ti{\pm 1},...,\ti{\pm 3}\}$ generically and locally give full dimensional tangent spaces to $\mathcal{L}_3$ and $\mathcal{L}_1$, respectively, at an intersection point $u$.  However, collectively they span a full dimensional tangent space to $\mathcal{X}$ since they are globally independent.  Hence, for any point in $u\in \mathcal{U}$, $T_u\mathcal{L}_1+T_u\mathcal{L}_3=T_u\mathcal{X}$.

Compounded with the fact that the leaves are properly transverse and trivially intersect, we conclude that for any $u\in \mathcal{U}$, $T_u\mathcal{X}=T_u\mathcal{L}_1\oplus T_u\mathcal{L}_3$; that is, the leaves are generically transverse.\end{proof}

We thus conclude that the tangent spaces $T_u\mathcal{X}$ are symplectic given by the product form.  This does not imply that $\mathcal{X}$ is complex symplectic since the form may not be closed.  We call two symplectic submanifolds symplectically dual if their tangent spaces are symplectic duals to each other with respect to this form.  Thus we have

\begin{corollary}
The symplectic leaves of $\mathcal{X}$ are (generically) symplectically dual.
\end{corollary}






\begin{commnt}
This sort of phenomena is not general.  For $\mathrm{SL}(2,\C)$ the leaves $\mathcal{L}(\Sigma_{3,0})$ and $\mathcal{L}(\Sigma_{1,1})$ are respectively dimension 0 and 2 and the variety $\X$ is dimension 3, so there is not transversality.  
\end{commnt}





\bibliography{bib}

\end{document}

%% file: oneholefrompants.arXiv.bbl
\providecommand{\bysame}{\leavevmode\hbox to3em{\hrulefill}\thinspace}
\providecommand{\MR}{\relax\ifhmode\unskip\space\fi MR }

\providecommand{\href}[2]{#2}